\documentclass[a4j,14pt]{article}  
\usepackage{geometry}
\usepackage[T1]{fontenc} 
\usepackage{times} 
\usepackage{amsmath,amssymb,bm, mathpazo}
\usepackage{amsfonts}
\usepackage{pifont}
\usepackage{mathrsfs}
\usepackage{graphicx}
\usepackage{amscd}
   
\makeatletter

\newtheorem{thm}{\theoremname}[section]

\newtheorem{df}[thm]{\definitionname}
\newtheorem{lemma}[thm]{\lemmaname}

\newtheorem{conj}[thm]{\conjecturename}

\newcommand{\theoremname}{Theorem}
\newcommand{\definitionname}{Definition}
\newcommand{\lemmaname}{Lemma}
\newcommand{\corollaryname}{Corollary}
\newcommand{\axiomname}{Axiom}
\newcommand{\propositionname}{Proposition}
\newcommand{\problemname}{Problem}
\newcommand{\examplename}{Example}
\newcommand{\remarkname}{Remark}

\newcommand{\conjecturename}{Conjecture}

\def\tedsymbol{\vcenter{\hbox{\vrule\@height.5em\@width.5em}}}
\def\ted{{\unskip\nobreak\hfil\penalty50
 \quad\hbox{}\nobreak\hfil \hbox{$\tedsymbol$}
 \parfillskip\z@ \finalhyphendemerits\z@\par}}

\makeatother

\title{On Previdi's delooping conjecture for $K$-theory}
\author{
\LARGE{Sho Saito}\footnote{Supported by JSPS Research Fellowships for Young Scientists.} 
}
\begin{document}
\date{\empty}
\maketitle
\nocite{Harada-01,Katura-05,Nakano-03,Takagi-96,Artin-74,Kuga-68,MP-64,Infeld-96,Borceux-01,Serre-95,Iyanaga-99,Iyanaga-02}
\begin{abstract}
We prove a modified version of Previdi's conjecture stating that the Waldhausen space ($K$-theory space) of an exact category is delooped by the Waldhausen space ($K$-theory space) of Beilinson's category of generalized Tate vector spaces. 
Our modified version states the delooping with non-connective $K$-theory spectra, almost including Previdi's original statement. 
As a consequence we obtain that the negative $K$-groups of an exact category are given by the $0$-th $K$-groups of the idempotent-completed iterated Beilinson categories, extending a theorem of Drinfeld on the first negative $K$-group.  
\end{abstract}
\section{Introduction} 
In his Ph.D. thesis \cite{previdi}, Previdi developes a categorical generalization of Kapranov's work \cite{kapranov} on dimensional and determinantal theories for Tate vector spaces over a field. 
His main results are formulated in terms of algebraic $K$-theory, and he observes certain relation between the $K$-groups $K_i(\mathcal{A})$, $i=0,1$, of an exact category $\mathcal{A}$ and $K_i(\displaystyle\lim_{\longleftrightarrow}\mathcal{A})$, $i=1,2$, of the exact category $\displaystyle\lim_{\longleftrightarrow}\mathcal{A}$ introduced by Beilinson \cite{beilinson}. 
He concludes the thesis with the following conjecture, which would include all the higher analogues of that relation. 
\begin{conj}[Previdi \cite{previdi}, 5.1.7] 
\label{theconjecture}
Write $S(\mathcal{A})$ for the geometric realization of the simplicial category $iS_{\bullet}(\mathcal{A})$ given by Waldhausen's $S_{\bullet}$-construction, the homotopy groups of whose loop space are the algebraic $K$-theory groups of the exact category $\mathcal{A}$. 
(See \cite{waldhausen}.) 
If $\mathcal{A}$ is partially abelian, i.e. if it and its opposite have pullbacks of admissible monomorphisms with common target, then $S(\mathcal{A})$ is delooped by $S(\displaystyle\lim_{\longleftrightarrow}\mathcal{A})$. 
\end{conj} 

In this article we prove the following modified version of the conjecture. 
\begin{thm}
\label{mainthm}
Let $\mathbb{K}(\mathcal{A})$ be the non-connective $K$-theory spectrum of the exact category $\mathcal{A}$, whose $i$-th homotopy group is the $i$-th $K$-group of $\mathcal{A}$ if $i>0$, the $0$-th $K$-group of the idempotent completion of $\mathcal{A}$ if $i=0$, and the $(-i)$-th negative $K$-group of $\mathcal{A}$ if $i<0$. 
(See \cite{schlichting2}.)  
Then there is a homotopy equivalence of spectra $\mathbb{K}(\mathcal{A})\xrightarrow[]{\sim}\Omega\mathbb{K}(\displaystyle\lim_{\longleftrightarrow}\mathcal{A})$. 
\end{thm} 
Note that no assumption on $\mathcal{A}$ is necessary. 
We also remark that Theorem \ref{mainthm} implies almost all of the essential part of Conjecture \ref{theconjecture}. 
Indeed, there results an isomorphism $K_i(\mathcal{A})\xrightarrow[]{\sim}K_{i+1}(\displaystyle\lim_{\longleftrightarrow}\mathcal{A})$ for any $\mathcal{A}$ and for every $i\geq 1$. 
If $\mathcal{A}$ is idempotent complete (this is the case for most of the typical examples such as the category $\mathcal{P}(R)$ of finitely generated projective modules over a ring $R$, and the category of 
vector bundles on a scheme, and any abelian category) this holds also for $i=0$. 
Theorem \ref{mainthm} moreover says that the $i$-th negative $K$-group $K_{-i}(\mathcal{A})$, $i>0$, is isomorphic to the $0$-th $K$-group of the idempotent completion of the $i$-times iterated Beilinson category $\displaystyle\lim_{\longleftrightarrow}{}^{i}\mathcal{A}$. 

\subsection*{Applications to the study of generalized Tate vector spaces}  
Previdi's work has its background in the study of generalized Tate vector spaces. 
Recall that a \textit{Tate vector space} over a discrete field $k$ is a topological $k$-vector space of the form $P\oplus Q^{\ast}$ where $P$ and $Q$ are discrete spaces and $(-)^{\ast}$ denotes the topological dual. 
The category of Tate vector spaces over $k$ is equivalent to the Beilinson category $\displaystyle\lim_{\longleftrightarrow}\operatorname{Vect}_0k$ of the exact category $\operatorname{Vect}_0k$ of finite-dimensional $k$-vector spaces. 
(See \cite{lcoec}, 7.4.) 

There are two directions of generalizing this notion, one of which due to Arkhipov-Kremnizer \cite{ak} et al. is the notion of an \textit{$n$-Tate vector space} as an object of the $n$-times iterated Beilinson category $\displaystyle\lim_{\longleftrightarrow}{}^n\operatorname{Vect}_0k$, $n\geq 1$. 
The other one due to Drinfeld \cite{drinfeld} generalizes the field $k$ to a commutative ring $R$ to get the notion of \textit{Tate $R$-module}. 
(We assume commutativity for simplicity, although Drinfeld's definition makes sense for non-commutative rings.) 
More precisely, he defines an \textit{elementary Tate $R$-module} to be a topological $R$-module of the form $P\oplus Q^{\ast}$, where $P$ and $Q$ are discrete projective $R$-modules, and a \textit{Tate $R$-module} to be a direct summand of an elementary Tate $R$-module. 
He shows that the notion of Tate $R$-module is local for the flat topology (\cite{drinfeld}, Theorem 3.3) and every Tate $R$-module is Nisnevich-locally elementary (\cite{drinfeld}, Theorem 3.4), and the first negative $K$-group $K_{-1}(R)$ of the ring $R$ is isomorphic to the $0$-th $K$-group of the exact category of Tate $R$-modules (\cite{drinfeld}, Theorem 3.6).   

The former one of these two generalizations is obtained purely formally by iterating the Beilinson construction, whereas the latter is based on nontrivial facts in ring theory.  
In fact, these two directions of generalization can be combined together. 
The equivalence of $\displaystyle\lim_{\longleftrightarrow}\operatorname{Vect}_0k$ to Tate $k$-vector spaces can be generalized to show that $\displaystyle\lim_{\longleftrightarrow}\mathcal{P}(R)$ is very similar to the category of elementary Tate $R$-modules, and hence the idempotent completion of $\displaystyle\lim_{\longleftrightarrow}\mathcal{P}(R)$ can be considered as a categorical substitute for Drinfeld's category of Tate $R$-modules. 
It is thus plausible to define an \textit{$n$-Tate $R$-module}, $n\geq 1$, to be an object of the idempotent completion of $\displaystyle\lim_{\longleftrightarrow}{}^n\mathcal{P}(R)$.  
Theorem \ref{mainthm} then can be regarded as a generalization of Drinfeld's theorem on $K_{-1}(R)$, as it says that the $n$-th negative $K$-group $K_{-n}(R)$ is isomorphic to the $0$-th $K$-group of $n$-Tate $R$-modules. 

Let us also discuss a consequence of Theorem \ref{mainthm} on $1$-Tate modules.  
Denote by $\mathcal{K}$ the sheaf of group-like $E_{\infty}$-spaces on the Nisnevich site of $\operatorname{Spec}R$, that sends an $R$-algebra $S$ to the space $\Omega^{\infty}\mathbb{K}(S)$. 
We describe how our Theorem \ref{mainthm} together with Drinfeld's theorem on the Nisnevich-local vanishing of $K_{-1}$ provides a purely formal way of associating to a $1$-Tate $R$-module $M$ a $\mathcal{K}$-torsor with a canonical action of the sheaf of groups of automorphisms of $M$.   
We note that this construction is essentially explained by Drinfeld (\cite{drinfeld}, section 5.5), who attributes it to Beilinson. 

Firstly, Theorem \ref{mainthm} shows that, in the $\infty$-topos of sheaves of spaces on the Nisnevich site of $\operatorname{Spec}R$, the sheaf $S\mapsto\Omega^{\infty}\mathbb{K}(\displaystyle\lim_{\longleftrightarrow}\mathcal{P}(S))$ is an object whose loop-space object is $\mathcal{K}$.  
It is obviously a pointed object. 
Drinfeld's theorem on the local vanishing of $K_{-1}$ in addition tells that this object is connected, i.e. $S\mapsto\Omega^{\infty}\mathbb{K}(\displaystyle\lim_{\longleftrightarrow}\mathcal{P}(S))$ is the classifying-space object for the $\infty$-group object $\mathcal{K}$. 
Then by general theory a $\mathcal{K}$-torsor corresponds to a map from the terminal object to the sheaf $S\mapsto\Omega^{\infty}\mathbb{K}(\displaystyle\lim_{\longleftrightarrow}\mathcal{P}(S))$, i.e. to a point of the space $\Omega^{\infty}\mathbb{K}(\displaystyle\lim_{\longleftrightarrow}\mathcal{P}(R))$.  
Thus the $1$-Tate $R$-module $M$, as an object of the idempotent completion of $\displaystyle\lim_{\longleftrightarrow}\mathcal{P}(R)$, defines such a torsor. 
The sheaf of groups of automorphisms of $M$ acts on it since, in general, for any idempotent complete exact category $\mathcal{A}$ and an object $A$ of $\mathcal{A}$ the classifying space of $\operatorname{Aut}_{\mathcal{A}}A$ admits a natural, canonical mapping to $\Omega S(\mathcal{A})=\Omega^{\infty}\mathbb{K}(\mathcal{A})$ which sends the base point to the point of $\Omega S(\mathcal{A})=\Omega^{\infty}\mathbb{K}(\mathcal{A})$ defined by the object $A$. 
(This is the composition of the map $B\operatorname{Aut}_{\mathcal{A}}A\to Bi\mathcal{A}$ with the first structure map $Bi\mathcal{A}\to\Omega S(\mathcal{A})$ of the connective algebraic $K$-theory spectrum of $\mathcal{A}$, where $i\mathcal{A}$ is the category of isomorphisms of $\mathcal{A}$, and $B$ indicates the classifying space of a category.)    

\subsection*{Organization and conventions} 
In section 2 we recall the definition and properties of the Beilinson category $\displaystyle\lim_{\longleftrightarrow}\mathcal{A}$ following Beilinson \cite{beilinson} and Previdi \cite{lcoec}. 
We recall the notions of ind- and pro-objects, introduce the categories $\operatorname{Ind}^a_{\mathbb{N}}\mathcal{A}$ and $\operatorname{Pro}^a_{\mathbb{N}}\mathcal{A}$, and discuss their relation to $\displaystyle\lim_{\longleftrightarrow}\mathcal{A}$. 
All statements in this section are either results of \cite{beilinson} and \cite{lcoec} or their immediate consequences. 

Section 3 begins with recalling Schlichting's results in \cite{schlichting1} which provides a powerful tool for constructing a homotopy fibration sequence of non-connective $K$-theory spectra. 
We prove Theorem \ref{mainthm} according to the following strategy: 
We construct, using Schlichting's method, two homotopy fibration sequences which fit into the commutative diagram 
\begin{displaymath}
\begin{CD}
\mathbb{K}(\mathcal{A})@>>>\mathbb{K}(\operatorname{Ind}_{\mathbb{N}}^a\mathcal{A})@>>>\mathbb{K}(\operatorname{Ind}_{\mathbb{N}}^a\mathcal{A}/\mathcal{A})\\
@VVV @VVV @VVV\\\
\mathbb{K}(\operatorname{Pro}_{\mathbb{N}}^a\mathcal{A})@>>>\mathbb{K}(\displaystyle \lim_{\longleftrightarrow}\mathcal{A})@>>>\mathbb{K}(\displaystyle \lim_{\longleftrightarrow}\mathcal{A}/\operatorname{Pro}_{\mathbb{N}}^a\mathcal{A}),   
\end{CD}
\end{displaymath}
as the horizontal sequences. 
We then go on to showing that the third vertical map is an equivalence, and that in the left hand square the upper right and lower left corners are contractible, so that the stated homotopy equivalence is obtained. 
(We remark that the upper horizontal homotopy fibration sequence and its consequence that $\mathbb{K}(\mathcal{A})$ is delooped by $\mathbb{K}(\operatorname{Ind}^a_{\mathbb{N}}\mathcal{A}/\mathcal{A})$ are Schlichting's results \cite{schlichting1}. 
Our delooping is a combination of his delooping with its dual.) 

We follow the notations adopted in \cite{previdi} and \cite{lcoec}. 
For instance we write $\operatorname{Ind}_{\mathbb{N}}^a\mathcal{A}$ for what is denoted by $\mathcal{F}\mathcal{A}$ in \cite{schlichting1}, and $\operatorname{Fun}^a(\Pi,\mathcal{A})$ instead of $\mathbf{A}^{\Pi}_a$ which is used in \cite{beilinson}.  
We write $\widetilde{\mathcal{A}}$ for the idempotent completion of $\mathcal{A}$. 
By saying a functor $\mathcal{A}\hookrightarrow\mathcal{U}$ between exact categories to be an \textit{embedding of exact categories}, we mean that it is a fully faithful exact functor whose essential image is closed under extensions in $\mathcal{U}$ and that a short sequence in $\mathcal{A}$ is exact if and only if its image in $\mathcal{U}$ is exact. 

\subsection*{Acknowledgement} 
First of all, I am greatly indebted to Marco Schlichting for many valuable ideas employed in the present article. 
In particular, the key idea of constructing the $s$-filtering localization sequence $\mathbb{K}(\operatorname{Pro}_{\mathbb{N}}^a\mathcal{A})\to\mathbb{K}(\displaystyle\lim_{\longleftrightarrow}\mathcal{A})\to\mathbb{K}(\displaystyle\lim_{\longleftrightarrow}\mathcal{A}/\operatorname{Pro}_{\mathbb{N}}^a\mathcal{A})$ and comparing it with the sequence $\mathbb{K}(\mathcal{A})\to\mathbb{K}(\operatorname{Ind}_{\mathbb{N}}^a\mathcal{A})\to\mathbb{K}(\operatorname{Ind}_{\mathbb{N}}^a\mathcal{A}/\mathcal{A})$ essentially belongs to him. 
I am truly grateful to Oliver Braeunling, Michael Groechnig and Jesse Wolfson for pointing out a crucial mistake in the previous version of the article. 
They moreover suggested how to fix the proof, letting me notice that any exact category admits a left $s$-filtering embedding into the category of admissible ind-objects and how the assumption of idempotent completeness can be removed in constructing the $s$-filtering localization sequences.  
I am also deeply grateful to Lars Hesselholt for commenting on early drafts of this work and for pointing out a misunderstanding of Previdi's result, and to Satoshi Mochizuki for answering many questions and suggesting valuable ideas in several times of helpful discussions. 
Finally, I thank Luigi Previdi for posing the nice conjecture and for answering my questions on its meaning. 

\section{Beilinson's category $\displaystyle \lim_{\longleftrightarrow}\mathcal{A}$} 
\subsection{Ind- and pro-objects in a category} 
We first recall some generalities on ind- and pro-objects. 
For any category $\mathcal{C}$, the category $\operatorname{Ind}\mathcal{C}$ (resp. $\operatorname{Pro}\mathcal{C}$) of \textit{ind-objects} (resp. \textit{pro-objects}) in $\mathcal{C}$ is defined to have as objects functors $\mathcal{X}:J\to\mathcal{C}$ with domain $J$ small and filtering (resp. $\mathcal{X}:I^{\operatorname{op}}\to\mathcal{C}$ with $I$ small and filtering). 
The ind-object $\mathcal{X}:J\to\mathcal{C}$ (resp. pro-object $\mathcal{X}:I^{\operatorname{op}}\to\mathcal{C}$) defines a functor $\mathcal{C}^{\operatorname{op}}\to(\text{sets})$, $C\mapsto\varinjlim_{j\in J}\operatorname{Hom}_{\mathcal{C}}(C,\mathcal{X}_j)$ (resp. $\mathcal{C}\to(\text{sets})$, $C\mapsto\varinjlim_{i\in I}\operatorname{Hom}_{\mathcal{C}}(\mathcal{X}_i,C)$). 
A morphism $\mathcal{X}\to \mathcal{Y}$ of ind-objects (resp. pro-objects) is a natural transformation between the functors $\mathcal{C}^{\operatorname{op}}\to(\text{sets})$ (resp. $\mathcal{C}\to(\text{sets})$) associated to $\mathcal{X}$ and $\mathcal{Y}$.   
Equivalently, the sets of morphisms of ind- and pro-objects can be defined to be the projective-inductive limits $\operatorname{Hom}_{\operatorname{Ind}\mathcal{C}}(\mathcal{X},\mathcal{Y})=\varprojlim_i\varinjlim_k\operatorname{Hom}_{\mathcal{C}}(\mathcal{X}_i,\mathcal{Y}_k)$ and $\operatorname{Hom}_{\operatorname{Pro}\mathcal{C}}(\mathcal{X},\mathcal{Y})=\varprojlim_l\varinjlim_j\operatorname{Hom}_{\mathcal{C}}(\mathcal{X}_j,\mathcal{Y}_l)$, respectively. 

If $\mathcal{X}$ and $\mathcal{Y}$ have a common index category, a natural transformation $\mathcal{X}\to\mathcal{Y}$ between the functors $\mathcal{X}$ and $\mathcal{Y}$ defines a map between the ind- or pro-objects $\mathcal{X}$ and $\mathcal{Y}$. 
Conversely, every map of ind- or pro-objects $\mathcal{X}\to\mathcal{Y}$ can be ``straightified'' to a natural transformation, in the sense that there is a commutative diagram in $\operatorname{Ind}\mathcal{C}$ or $\operatorname{Pro}\mathcal{C}$ 
\begin{displaymath}
\begin{CD}
\mathcal{X}@>>>\mathcal{Y}\\
@V{\sim}VV@V{\sim}VV\\
\widetilde{\mathcal{X}}@>>>\widetilde{\mathcal{Y}}
\end{CD}
\end{displaymath} 
with the vertical maps isomorphisms, $\widetilde{\mathcal{X}}$ and $\widetilde{\mathcal{Y}}$ having a common index category, and $\widetilde{\mathcal{X}}\to\widetilde{\mathcal{Y}}$ coming from a natural transformation.  
(See \cite{am}, Appendix, for details.) 

If $\mathcal{C}$ is an exact category, the categories $\operatorname{Ind}\mathcal{C}$ and $\operatorname{Pro}\mathcal{C}$ possess exact structures. 
A pair of composable morphisms in $\operatorname{Ind}\mathcal{C}$ or $\operatorname{Pro}\mathcal{C}$ is a short exact sequence if it can be straightified to a sequence of natural transformations which is level-wise exact in $\mathcal{C}$ (\cite{lcoec}, 4.15, 4.16).   
In this article we are mainly concerned with the full subcategories $\operatorname{Ind}^a\mathcal{C}$ and $\operatorname{Pro}^a\mathcal{C}$ of \textit{admissible ind-} and \textit{pro-objects} introduced by Previdi \cite{lcoec}, 5.6: 
An ind-object $\mathcal{X}:J\to\mathcal{C}$ (resp. pro-object $\mathcal{X}:I^{\operatorname{op}}\to\mathcal{C}$) is \textit{admissible} if for every map $j\to j^{\prime}$ in $J$ (resp. $i\to i^{\prime}$ in $I$) the morphism $X_j\hookrightarrow X_{j^{\prime}}$ is an admissible monomorphism in $\mathcal{C}$ (resp. $X_i\twoheadleftarrow X_{i^{\prime}}$ an admissible epimorphism).  
These subcategories are extension-closed in the exact categories $\operatorname{Ind}\mathcal{C}$ and $\operatorname{Pro}\mathcal{C}$, respectively, so that they have induced exact structures. 
Since an object $C$ of $\mathcal{C}$ can be considered as an admissible ind- or pro-object which is indexed by whatever small and filtering category and takes the constant value $C$, there are embeddings of exact categories $\mathcal{C}\hookrightarrow\operatorname{Ind}^a\mathcal{C}$ and $\mathcal{C}\hookrightarrow\operatorname{Pro}^a\mathcal{C}$.  

We write $\operatorname{Ind}^a_{\mathbb{N}}\mathcal{C}$ and $\operatorname{Pro}^a_{\mathbb{N}}\mathcal{C}$ for the full, extension-closed subcategories of $\operatorname{Ind}^a\mathcal{C}$ and $\operatorname{Pro}^a\mathcal{C}$ consisting of admissible ind- and pro-objects, respectively, indexed by the filtering category of natural numbers. 
(There is precisely one morphism $j\to k$ if $j\leq k\in\mathbb{N}$.) 
The object $C$ of $\mathcal{C}$ defines an object $C=C=C=\cdots$ in $\operatorname{Ind}^a_{\mathbb{N}}\mathcal{C}$ or $\operatorname{Pro}^a_{\mathbb{N}}\mathcal{C}$. 
Note that the resulting embedding $\mathcal{C}\hookrightarrow\operatorname{Ind}^a_{\mathbb{N}}\mathcal{C}\hookrightarrow\operatorname{Ind}^a\mathcal{C}$ (resp. $\mathcal{C}\hookrightarrow\operatorname{Pro}^a_{\mathbb{N}}\mathcal{C}\hookrightarrow\operatorname{Pro}^a\mathcal{C}$) is naturally isomorphic to the embedding $\mathcal{C}\hookrightarrow\operatorname{Ind}^a\mathcal{C}$ (resp. $\mathcal{C}\hookrightarrow\operatorname{Pro}^a\mathcal{C}$) mentioned above. 

\subsection{Definition of $\displaystyle \lim_{\longleftrightarrow}\mathcal{A}$} 
Let $\mathcal{A}$ be an exact category. 
We write $\Pi$ for the ordered set $\{(i,j)\in\mathbb{Z}\times\mathbb{Z}\mid i\leq j\}$, where $(i,j)\leq(i^{\prime},j^{\prime})$ if $i\leq i^{\prime}$ and $j\leq j^{\prime}$. 
A functor $X:\Pi\to\mathcal{A}$, where $\Pi$ is viewed as a filtered category, is {\em admissible} if for every triple $i\leq j\leq k$, the sequence $X_{i,j}\hookrightarrow X_{i,k}\twoheadrightarrow X_{j,k}$ is a short exact sequence in $\mathcal{A}$. 
We denote by $\operatorname{Fun}^a(\Pi,\mathcal{A})$ the exact category of admissible functors $X:\Pi\to\mathcal{A}$ and natural transformations, where a short sequence $X\to Y\to Z$ of natural transformations of admissible functors is a short exact sequence in $\operatorname{Fun}^a(\Pi,\mathcal{A})$ if $X_{i,j}\hookrightarrow Y_{i,j}\twoheadrightarrow Z_{i,j}$ is a short exact sequence in $\mathcal{A}$ for every $i\leq j$. 
A bicofinal map $\phi:\mathbb{Z}\to\mathbb{Z}$ ($\phi$ is said to be \textit{bicofinal} if it is nondecreasing and satisfies $\lim_{i\to\pm\infty}\phi(i)=\pm\infty$) induces a cofinal functor $\widetilde{\phi}:\Pi\to\Pi$, $(i,j)\mapsto(\phi(i),\phi(j))$.  
If $\phi$ and $\psi:\mathbb{Z}\to\mathbb{Z}$ are bicofinal maps such that $\phi(i)\leq\psi(i)$ for all $i$, and if $X:\Pi\to\mathcal{A}$ is an admissible functor, there is a natural transformation $u_{X,\phi,\psi}:X\circ\widetilde{\phi}\to X\circ\widetilde{\psi}$. 
\begin{df}[Beilinson \cite{beilinson}, A.3] 
The category $\displaystyle \lim_{\longleftrightarrow}\mathcal{A}$ is defined to be the localization of $\operatorname{Fun}^a(\Pi,\mathcal{A})$ by the morphisms $u_{X,\phi,\psi}$, where $X\in\operatorname{ob}\operatorname{Fun}^a(\Pi,\mathcal{A})$, and $\phi\leq\psi:\mathbb{Z}\to\mathbb{Z}$ are bicofinal. 
\end{df}  

If $X:\Pi\to\mathcal{A}$ is an admissible functor, we have for each $j\in\mathbb{Z}$ an admissible pro-object $X_{\bullet,j}:\{i\in\mathbb{Z}\mid i\leq j\}\to\mathcal{A}$, $i\mapsto X_{i,j}$, in $\mathcal{A}$. 
We get in turn an admissible ind-object $\mathbb{Z}\to\operatorname{Pro}^a\mathcal{A}$, $j\mapsto X_{\bullet,j}$, in $\operatorname{Pro}^a\mathcal{A}$.    
Thus the admissible functor $X$ can be viewed as an object of the iterated $\operatorname{Ind}$-$\operatorname{Pro}$ category $\operatorname{Ind}^a\operatorname{Pro}^a\mathcal{A}$. 
If $\phi\leq\psi:\mathbb{Z}\to\mathbb{Z}$ are bicofinal, the map $u_{X,\phi,\psi}$ defines an isomorphism between the ind-pro-objects $X\circ\widetilde{\phi}$ and $X\circ\widetilde{\psi}$. 
We get a functor $\displaystyle \lim_{\longleftrightarrow}\mathcal{A}\to\operatorname{Ind}^a\operatorname{Pro}^a\mathcal{A}$.     
In view of the following theorem, we regard $\displaystyle \lim_{\longleftrightarrow}\mathcal{A}$ as an exact subcategory of $\operatorname{Ind}^a\operatorname{Pro}^a\mathcal{A}$.   
\begin{thm}[Previdi \cite{lcoec}, 5.8, 6.1] 
The functor $\displaystyle \lim_{\longleftrightarrow}\mathcal{A}\to\operatorname{Ind}^a\operatorname{Pro}^a\mathcal{A}$ is fully faithful. 
Moreover, the image is closed under extensions in $\operatorname{Ind}^a\operatorname{Pro}^a\mathcal{A}$. 
In particular, $\displaystyle\lim_{\longleftrightarrow}\mathcal{A}$ has an exact structure where a sequence in $\displaystyle \lim_{\longleftrightarrow}\mathcal{A}$ is exact if and only if its image in $\operatorname{Ind}^a\operatorname{Pro}^a\mathcal{A}$ is exact. 
\end{thm}   
By \cite{lcoec}, 6.3, there is an embedding $\operatorname{Ind}^a_{\mathbb{N}}\mathcal{A}\hookrightarrow\displaystyle \lim_{\longleftrightarrow}\mathcal{A}$ (resp. $\operatorname{Pro}^a_{\mathbb{N}}\mathcal{A}\hookrightarrow\displaystyle \lim_{\longleftrightarrow}\mathcal{A}$) of exact categories that sends $X_1\hookrightarrow X_2\hookrightarrow X_3\hookrightarrow\cdots\in\operatorname{ob}\operatorname{Ind}^a_{\mathbb{N}}\mathcal{A}$ (resp. $X_1\twoheadleftarrow X_2\twoheadleftarrow X_3\twoheadleftarrow\cdots\in\operatorname{Pro}^a_{\mathbb{N}}\mathcal{A}$) to the object in $\displaystyle \lim_{\longleftrightarrow}\mathcal{A}$ determined by $X_{i,j}=X_{0,j}=X_j$ for $i\leq 0<j$ (resp. $X_{i,j}=X_{i,1}=X_{-i+1}$ for $i\leq 0<j$).  

We refer to \cite{lcoec} for detailed discussions on ind/pro-objects in exact categories.  

\section{Proof of Theorem \ref{mainthm}} 
We prove the theorem using the $s$-filtering localization sequence constructed by Schlichting \cite{schlichting1}. 

Let $\mathcal{A}\hookrightarrow\mathcal{U}$ be an embedding of exact categories.  
Following Schlichting \cite{schlichting1}, we define a map in $\mathcal{U}$ to be a {\em weak isomorphism} with respect to $\mathcal{A}\hookrightarrow\mathcal{U}$ if it is either an admissible monomorphism that admits a cokernel in the essential image of $\mathcal{A}\hookrightarrow\mathcal{U}$ or an admissible epimorphism that admits a kernel in the essential image of $\mathcal{A}\hookrightarrow\mathcal{U}$. 
In particular, for every $A\in\operatorname{ob}\mathcal{A}$ the maps $0\to A$ and $A\to 0$ are weak isomorphisms.  
The localization of $\mathcal{U}$ by weak isomorphisms with respect to $\mathcal{A}$ is denoted by $\mathcal{U}/\mathcal{A}$.  
Recall, from \cite{schlichting1}, that the embedding $\mathcal{A}\hookrightarrow\mathcal{U}$ of exact categories is a {\em left $s$-filtering} if the following conditions are satisfied. 
\begin{itemize}
\item[(1)] If $A\twoheadrightarrow U$ is an admissible epimorphism in $\mathcal{U}$ with $A\in\operatorname{ob}\mathcal{A}$, then $U\in\operatorname{ob}\mathcal{A}$. 
\item[(2)] If $U\hookrightarrow A$ is an admissible monomorphism in $\mathcal{U}$ with $A\in\operatorname{ob}\mathcal{A}$, then $U\in\operatorname{ob}\mathcal{A}$. 
\item[(3)] Every map $A\to U$ in $\mathcal{U}$ with $A\in\operatorname{ob}\mathcal{A}$ factors through an object $B\in\operatorname{ob}\mathcal{A}$ such that $B\hookrightarrow U$ is an admissible monomorphism in $\mathcal{U}$. 
\item[(4)] If $U\twoheadrightarrow A$ is an admissible epimorphism in $\mathcal{U}$ with $A\in\operatorname{ob}\mathcal{A}$, then there is an admissible monomorphism $B\hookrightarrow U$ with $B\in\operatorname{ob}\mathcal{A}$ such that the composition $B\twoheadrightarrow A$ is an admissible epimorphism in $\mathcal{A}$. 
\end{itemize} 
(Here $\operatorname{ob}\mathcal{A}$ denotes by slight abuse of notation the collection of objects of $\mathcal{U}$ contained in the essential image of $\mathcal{A}\hookrightarrow\mathcal{U}$.) 
A {\em right $s$-filtering} embedding is defined by dualizing the conditions above. 

The following theorem due to Schlichting \cite{schlichting1} we use as the main technical tool for the proof. 
\begin{thm}[Schlichting \cite{schlichting1}, 1.16, 1.20, 2.10] 
\label{schlichtingfibration} 
If $\mathcal{A}\hookrightarrow\mathcal{U}$ is left or right $s$-filtering, then the localization $\mathcal{U}/\mathcal{A}$ has an exact structure where a short sequence is exact if and only if it is isomorphic to the image of a short exact sequence in $\mathcal{U}$. 
Moreover, the sequence of exact categories $\mathcal{A}\to\mathcal{U}\to\mathcal{U}/\mathcal{A}$ induces a homotopy fibration $\mathbb{K}(\mathcal{A})\to\mathbb{K}(\mathcal{U})\to\mathbb{K}(\mathcal{U}/\mathcal{A})$ of non-connective $K$-theory spectra. 
\end{thm} 
\textit{Remark}. 
The single statement of Theorem 2.10 of \cite{schlichting1} assumes the idempotent completeness of $\mathcal{A}$ for this $K$-theory sequence to be a homotopy fibration.  
But the theorem holds for general $\mathcal{A}$ in view of Lemma 1.20 of loc. cit., which assures whenever $\mathcal{A}\hookrightarrow\mathcal{U}$ is left or right $s$-filtering the existence of an extension-closed full subcategory $\widetilde{\mathcal{U}}^{\mathcal{A}}$ of $\widetilde{\mathcal{U}}$ such that $\mathcal{U}$ is cofinally contained in $\widetilde{\mathcal{U}}^{\mathcal{A}}$, the induced embedding $\widetilde{\mathcal{A}}\hookrightarrow\widetilde{\mathcal{U}}$ factors through a left or right $s$-filtering embedding $\widetilde{\mathcal{A}}\hookrightarrow\widetilde{\mathcal{U}}^{\mathcal{A}}$, and $\mathcal{U}/\mathcal{A}\xrightarrow[]{\sim}\widetilde{\mathcal{U}}^{\mathcal{A}}/\widetilde{\mathcal{A}}$ is an equivalence of exact categories. 
The homotopy fibration sequence $\mathbb{K}(\widetilde{\mathcal{A}})\to\mathbb{K}(\widetilde{\mathcal{U}})\to\mathbb{K}(\widetilde{\mathcal{U}}^{\mathcal{A}}/\widetilde{\mathcal{A}})$ is equivalent to the sequence $\mathbb{K}(\mathcal{A})\to\mathbb{K}(\mathcal{U})\to\mathbb{K}(\mathcal{U}/\mathcal{A})$ since a cofinal embedding of exact categories induces an equivalence of non-connective $K$-theory spectra.  

\begin{lemma} 
\label{sfilt}
For any exact category $\mathcal{A}$, the embedding $\mathcal{A}\hookrightarrow\operatorname{Ind}^a\mathcal{A}$ is left $s$-filtering. 
\end{lemma}
\begin{proof}
We start by checking condition (3) of left $s$-filtering. 
Let $X$ be an object of $\mathcal{A}$ and $Y$ an admissible ind-object in $\mathcal{A}$ indexed by a small filtering category $J$. 
A morphism $f:X\to Y$ in $\operatorname{Ind}^a\mathcal{A}$ is an element of $\varinjlim_{j\in J}\operatorname{Hom}_{\mathcal{A}}(X, Y_j)$, i.e. represented as the class of a map $f_j:X\to Y_j$ in $\mathcal{A}$ for some $j\in J$.  
The canonical map $Y_j\hookrightarrow Y$ is an admissible monomorphism because the diagram $j/J\to\mathcal{A}$, $i\mapsto Y_i/Y_j$ serves as its cockerel, where $j/J$ is the under-category of $j$. 
We get a desired factorization $f:X\xrightarrow[]{f_j}Y_j\hookrightarrow Y$.  

Condition (1) follows from (3). 
Indeed, an admissible epimorphism $X\twoheadrightarrow Y$ with $X$ in $\mathcal{A}$ factors through some $Z$ in $\mathcal{A}$ such that $Z\hookrightarrow Y$ is an admissible monomorphism. 
The composition $X\twoheadrightarrow Y\twoheadrightarrow Y/Z$ is $0$, but since this composition is also an admissible epimorphism, $Y/Z$ must be $0$. 
This forces $Y$ to be essentially constant. 

To prove (4), let $Y\twoheadrightarrow X$ be an admissible epimorphism in $\operatorname{Ind}^a\mathcal{A}$ with $X$ in $\mathcal{A}$, whose kernel we denote by $Z$. 
The short exact sequence $0\to Z\hookrightarrow Y\twoheadrightarrow X\to0$ is isomorphic to a straight exact sequence $0\to Z^{\prime}\hookrightarrow Y^{\prime}\twoheadrightarrow X^{\prime}\to0$, where $Z^{\prime}$, $Y^{\prime}$, and $X^{\prime}$ are all indexed by an identical small filtering category $I$ and respectively isomorphic to $Z$, $Y$, and $X$. 
The isomorphism $X^{\prime}\xrightarrow[]{\sim}X$ is a compatible collection of morphisms $g_i:X^{\prime}_i\to X$ in $\mathcal{A}$, $i\in I$, such that there is a morphism $h:X\to X^{\prime}_{i_0}$ for some $i_0\in I$ such that $g_{i_0}\circ h=\operatorname{id}_X$ and $h\circ g_{i_0}$ is equivalent to $\operatorname{id}_{X^{\prime}_{i_0}}$ in $\varinjlim_{i\in I}\operatorname{Hom}_{\mathcal{A}}(X^{\prime}_{i_0}, X^{\prime}_i)$. 
Since $X^{\prime}$ is an admissible ind-object this implies $h\circ g_{i_0}=\operatorname{id}_{X^{\prime}_{i_0}}$, i.e. $g_{i_0}$ is an isomorphism. 
(Note also that $g_i$ are isomorphisms for all $i\in i_0/I$.) 
The map $Y^{\prime}_{i_0}\hookrightarrow Y^{\prime}\xrightarrow[]{\sim}Y$ is an admissible monomorphism as noted above, and its composition with $Y\twoheadrightarrow X$ equals the composition $Y^{\prime}_{i_0}\twoheadrightarrow X^{\prime}_{i_0}\xrightarrow[g_{i_0}]{\sim}X$, which is an admissible epimorphism in $\mathcal{A}$. 

Finally, if $Y\hookrightarrow X$ is an admissible monomorphism with $X$ in $\mathcal{A}$, its cokernel $Z$ is in $\mathcal{A}$ by (1). 
Let $0\to Y^{\prime}\hookrightarrow X^{\prime}\twoheadrightarrow Z^{\prime}\to0$ be a straightification of the exact sequence $0\to Y\hookrightarrow X\twoheadrightarrow Z\to0$, whose common indices we denote by $I$. 
Then an argument similar to above shows that there is an $i_0\in I$ such that for every $i\in i_0/I$, $X^{\prime}_i$ and $Z^{\prime}_i$ are isomorphic to $X$ and $Z$, respectively. 
It follows that $Y^{\prime}_i$ is essentially constant above $i_0$, and we conclude that $Y$ is contained in the essential image of $\mathcal{A}$, verifying condition (2).  
\end{proof}

We remark that, given a composable pair of embeddings of exact categories $\mathcal{A}\hookrightarrow\mathcal{V}$ and $\mathcal{V}\hookrightarrow\mathcal{U}$, if their composition is naturally isomorphic to a left $s$-filtering embedding $\mathcal{A}\hookrightarrow\mathcal{U}$ then $\mathcal{A}\hookrightarrow\mathcal{V}$ is also left $s$-filtering. 
This in particular implies that the embeddings $\mathcal{A}\hookrightarrow\operatorname{Ind}^a_{\mathbb{N}}\mathcal{A}$ and $\operatorname{Pro}^a_{\mathbb{N}}\mathcal{A}\hookrightarrow\displaystyle\lim_{\longleftrightarrow}\mathcal{A}$ are left $s$-filtering. 
Hence by Theorem \ref{schlichtingfibration} we get two homotopy fibration sequences of non-connective $K$-theory spectra $\mathbb{K}(\mathcal{A})\to\mathbb{K}(\operatorname{Ind}^a_{\mathbb{N}}\mathcal{A})\to\mathbb{K}(\operatorname{Ind}^a_{\mathbb{N}}\mathcal{A}/\mathcal{A})$ and $\mathbb{K}(\operatorname{Pro}^a_{\mathbb{N}}\mathcal{A})\to\mathbb{K}(\displaystyle\lim_{\longleftrightarrow}\mathcal{A})\to\mathbb{K}(\displaystyle\lim_{\longleftrightarrow}\mathcal{A}/\operatorname{Pro}^a_{\mathbb{N}}\mathcal{A})$. 
We compare these sequences to deduce Theorem \ref{mainthm}. 
\begin{lemma} 
\label{equivalence} 
There is an equivalence $\displaystyle \operatorname{Ind}^a_{\mathbb{N}}\mathcal{A}/\mathcal{A}\xrightarrow[]{\sim}\lim_{\longleftrightarrow}\mathcal{A}/\operatorname{Pro}^a_{\mathbb{N}}\mathcal{A}$. 
\end{lemma}
\begin{proof}
We have a commutative diagram 
\begin{displaymath}
\begin{CD}
\mathcal{A}@>>>\operatorname{Ind}^a_{\mathbb{N}}\mathcal{A}\\ 
@VVV@VVV\\ 
\operatorname{Pro}^a_{\mathbb{N}}\mathcal{A}@>>>\displaystyle \lim_{\longleftrightarrow}\mathcal{A},  
\end{CD}
\end{displaymath}
whence there results a functor $\displaystyle F:\operatorname{Ind}^a_{\mathbb{N}}\mathcal{A}/\mathcal{A}\to\lim_{\longleftrightarrow}\mathcal{A}/\operatorname{Pro}^a_{\mathbb{N}}\mathcal{A}$. 

To construct a quasi inverse, we start by noticing that the functor $\operatorname{Fun}^a(\Pi,\mathcal{A})\to\operatorname{Ind}^a_{\mathbb{N}}\mathcal{A}$, $(X_{i,j})_{i\leq j}\mapsto X_{0,1}\hookrightarrow X_{0,2}\hookrightarrow\cdots$, induces a functor $\widetilde{G}:\displaystyle \lim_{\longleftrightarrow}\mathcal{A}\to\operatorname{Ind}^a_{\mathbb{N}}\mathcal{A}/\mathcal{A}$.  
Indeed, if $\phi\leq\psi:\mathbb{Z}\to\mathbb{Z}$ are bicofinal, the map $u_{X,\phi,\psi}:X\circ\widetilde{\phi}\to X\circ\widetilde{\psi}$ in $\operatorname{Fun}^a(\Pi,\mathcal{A})$ is sent to the map $X_{\phi(0),\phi(\bullet)}\to X_{\psi(0),\psi(\bullet)}$, which factors as $X_{\phi(0),\phi(\bullet)}\hookrightarrow X_{\phi(0),\psi(\bullet)}\twoheadrightarrow X_{\psi(0),\psi(\bullet)}$. 
The map $X_{\phi(0),\phi(\bullet)}\hookrightarrow X_{\phi(0),\psi(\bullet)}$ is an isomorphism in $\operatorname{Ind}_{\mathbb{N}}^{a}\mathcal{A}$ since it consists of natural isomorphisms $\varinjlim_{j}\operatorname{Hom}_{\mathcal{A}}(A,X_{\phi(0),\phi(j)})\xrightarrow[]{\sim}\varinjlim_{j}\operatorname{Hom}_{\mathcal{A}}(A,X_{\phi(0),\psi(j)})$, $A\in\operatorname{ob}\mathcal{A}$, as $\phi$ and $\psi$ are bicofinal.    
We also see that $X_{\phi(0),\psi(\bullet)}\twoheadrightarrow X_{\psi(0),\psi(\bullet)}$ is a weak isomorphism in $\operatorname{Ind}_{\mathbb{N}}^a\mathcal{A}$ with respect to $\mathcal{A}$, since it has the constant kernel $X_{\phi(0),\psi(0)}=X_{\phi(0),\psi(0)}=\cdots$.       
The functor $\widetilde{G}$ thus defined takes weak isomorphisms in $\displaystyle \lim_{\longleftrightarrow}\mathcal{A}$ with respect to $\operatorname{Pro}^a_{\mathbb{N}}\mathcal{A}$ to weak isomorphisms in $\operatorname{Ind}^a_{\mathbb{N}}\mathcal{A}$ with respect to $\mathcal{A}$, since if $\displaystyle X\in\operatorname{ob}\lim_{\longleftrightarrow}\mathcal{A}$ is in the image of $\operatorname{Pro}^a_{\mathbb{N}}\mathcal{A}$ then its $0$-th row is constant $X_{0,1}=X_{0,1}=\cdots$, i.e. $\widetilde{G}(X)$ is in the image of $\mathcal{A}$. 
Hence $\widetilde{G}$ factors through a functor $\displaystyle G:\lim_{\longleftrightarrow}\mathcal{A}/\operatorname{Pro}^a_{\mathbb{N}}\mathcal{A}\to\operatorname{Ind}^a_{\mathbb{N}}\mathcal{A}/\mathcal{A}$. 

We have $G\circ F=\operatorname{id}_{\operatorname{Ind}^a_{\mathbb{N}}\mathcal{A}/\mathcal{A}}$ by definition. 
On the other hand, if $X=\displaystyle (X_{i,j})_{i\leq j}\in\operatorname{ob}\lim_{\longleftrightarrow}\mathcal{A}$, then $F\circ G(X)$ is the object $\widetilde{X}$ of $\displaystyle \lim_{\longleftrightarrow}\mathcal{A}$ determined by $\widetilde{X}_{i,j}=\widetilde{X}_{0,j}=X_{0,j}$, ($i\leq 0<j$). 
Define an admissible epimorphism $f_X:X\twoheadrightarrow\widetilde{X}$ in $\operatorname{Fun}^a(\Pi,\mathcal{A})$ (hence in $\displaystyle \lim_{\longleftrightarrow}\mathcal{A}$) by 
\begin{displaymath}
(f_X)_{i,j}=
\begin{cases}
X_{i,j}=X_{i,j}&(0\leq i\leq j)\\
X_{i,j}\twoheadrightarrow X_{0,j}&(i\leq 0<j)\\
X_{i,j}\twoheadrightarrow 0&(i\leq j\leq 0).   
\end{cases}
\end{displaymath}  
The kernel coincides with the image of $0\twoheadleftarrow X_{-1,0}\twoheadleftarrow X_{-2,0}\twoheadleftarrow X_{-3,0}\twoheadleftarrow\cdots\in\operatorname{ob}\operatorname{Pro}_{\mathbb{N}}^a\mathcal{A}$ in $\displaystyle \lim_{\longleftrightarrow}\mathcal{A}$.  
Hence $f_X$ is a weak isomorphism in $\displaystyle \lim_{\longleftrightarrow}\mathcal{A}$ with respect to $\operatorname{Pro}^a_{\mathbb{N}}\mathcal{A}$. 
Thus we get an isomorphism $f:\operatorname{id}_{\displaystyle \lim_{\longleftrightarrow}\mathcal{A}/\operatorname{Pro}^a_{\mathbb{N}}\mathcal{A}}\xrightarrow[]{\sim}F\circ G$, to conclude that $G$ is a quasi inverse to $F$. 
\end{proof}

This means that in the commutative diagram of non-connective $K$-theory spectra  
\begin{displaymath}
\begin{CD}
\mathbb{K}(\mathcal{A})@>>>\mathbb{K}(\operatorname{Ind}_{\mathbb{N}}^a\mathcal{A})@>>>\mathbb{K}(\operatorname{Ind}_{\mathbb{N}}^a\mathcal{A}/\mathcal{A})\\
@VVV @VVV @VVV\\\
\mathbb{K}(\operatorname{Pro}_{\mathbb{N}}^a\mathcal{A})@>>>\mathbb{K}(\displaystyle \lim_{\longleftrightarrow}\mathcal{A})@>>>\mathbb{K}(\displaystyle \lim_{\longleftrightarrow}\mathcal{A}/\operatorname{Pro}_{\mathbb{N}}^a\mathcal{A}), 
\end{CD}
\end{displaymath}
the third vertical map is an equivalence. 
Since the two horizontal sequences are homotopy fibrations, it follows that the square 
\begin{displaymath}
\begin{CD}
\mathbb{K}(\mathcal{A})@>>>\mathbb{K}(\operatorname{Ind}_{\mathbb{N}}^a\mathcal{A})\\
@VVV @VVV\\
\mathbb{K}(\operatorname{Pro}_{\mathbb{N}}^a\mathcal{A})@>>>\mathbb{K}(\displaystyle \lim_{\longleftrightarrow}\mathcal{A})
\end{CD}
\end{displaymath}
is homotopy cartesian, i.e. $\mathbb{K}(\mathcal{A})\xrightarrow[]{\sim}\operatorname{holim}(\mathbb{K}(\operatorname{Pro}^a_{\mathbb{N}}\mathcal{A})\rightarrow\mathbb{K}(\displaystyle\lim_{\longleftrightarrow}\mathcal{A})\leftarrow\mathbb{K}(\operatorname{Ind}^a_{\mathbb{N}}\mathcal{A}))$ is an equivalence.  
We finally note:  
\begin{lemma}
\label{contractible} 
There are canonical contractions for the non-connective $K$-theory spectra $\mathbb{K}(\operatorname{Ind}_{\mathbb{N}}^a\mathcal{A})$ and $\mathbb{K}(\operatorname{Pro}_{\mathbb{N}}^a\mathcal{A})$. 
\end{lemma}
\begin{proof}
The contraction for $\mathbb{K}(\operatorname{Ind}^a_{\mathbb{N}}\mathcal{A})$ comes from the canonical flasque structure on $\operatorname{Ind}^a_{\mathbb{N}}\mathcal{A}$ (i.e. an endo-functor whose direct sum with the identity functor is naturally isomorphic to itself) given as follows. 
Let $X=(X_i)_{i\geq 1}\in\operatorname{ob}\operatorname{Ind}^a_{\mathbb{N}}\mathcal{A}$ be an $\mathbb{N}$-indexed admissible ind-object in $\mathcal{A}$ whose structure maps we denote by $\rho=\rho_{i,i^{\prime}}:X_i\hookrightarrow X_{i^{\prime}}$. 
Write $T(X)\in\operatorname{ob}\operatorname{Ind}^a_{\mathbb{N}}\mathcal{A}$ for the admissible ind-object $0\xrightarrow[]{}X_1\xrightarrow[]{(\rho,0)}
X_2\oplus X_1\xrightarrow[]{(\rho\oplus\rho,0)}X_3\oplus X_2\oplus X_1\xrightarrow[]{(\rho\oplus\rho\oplus\rho,0)}\cdots$. 
A morphism $f\in\operatorname{Hom}_{\operatorname{Ind}^a_{\mathbb{N}}\mathcal{A}}(Y,X)=\varprojlim_i\varinjlim_j\operatorname{Hom}_{\mathcal{A}}(Y_i,X_j)$ with $i$-th component represented by $f_i:Y_i\to X_{j(i)}$ defines a morphism $T(f):T(Y)\to T(X)$  whose $i$-th component is the class of the composition $Y_{i-1}\oplus\cdots\oplus Y_1\xrightarrow[]{f_{i-1}\oplus\cdots\oplus f_1}X_{j(i-1)}\oplus\cdots\oplus X_{j(1)}\xrightarrow[]{\rho\oplus\cdots\oplus\rho}X_{k+i-1}\oplus\cdots\oplus X_{k+1}\hookrightarrow T(X)_{k+i}$, where $k$ is chosen to be sufficiently large. 
The endo-functor $T$ thus defined is a flasque structure on $\operatorname{Ind}^a_{\mathbb{N}}\mathcal{A}$ since $(X\oplus T(X))_i\xrightarrow[]{=}T(X)_{i+1}$ give a natural isomorphism of ind-objects. 

The contraction for $\mathbb{K}(\operatorname{Pro}^a_{\mathbb{N}}\mathcal{A})$ follows from the contraction for $\mathbb{K}(\operatorname{Ind}^a_{\mathbb{N}}(-))$ via the identification $\operatorname{Pro}^a_{\mathbb{N}}\mathcal{A}=(\operatorname{Ind}^a_{\mathbb{N}}\mathcal{A}^{\operatorname{op}})^{\operatorname{op}}$ and the general equivalence $\mathbb{K}(\mathcal{B}^{\operatorname{op}})\xrightarrow[]{\sim}\mathbb{K}(\mathcal{B})$. 
\end{proof}

We now obtain the desired homotopy equivalence 
$\mathbb{K}(\mathcal{A})=\operatorname{holim}(\mathbb{K}(\operatorname{Pro}_{\mathbb{N}}^a\mathcal{A})\to\mathbb{K}(\displaystyle\lim_{\longleftrightarrow}\mathcal{A})\leftarrow\mathbb{K}(\operatorname{Ind}_{\mathbb{N}}^a\mathcal{A}))\stackrel{\sim}{\to}\operatorname{holim}(\ast\to\mathbb{K}(\displaystyle\lim_{\longleftrightarrow}\mathcal{A})\leftarrow\ast)=\Omega\mathbb{K}(\displaystyle\lim_{\longleftrightarrow}\mathcal{A})$, and the proof of Theorem \ref{mainthm} is complete. 
 
\begin{flushleft} 
Sho Saito\\
Graduate School of Mathematics\\
Nagoya University\\
Nagoya, Japan\\
\texttt{sho.saito@math.nagoya-u.ac.jp}
\end{flushleft}
\end{document}